 \documentclass[11pt]{amsart}
\usepackage{epsfig,graphicx}

\newcommand{\mcg}{\Gamma _{g,r}^n}

\newcommand{\M}{{\Gamma }}

\newtheorem{theorem}{Theorem}[section]
\newtheorem{thm}[theorem]{Theorem}

\newtheorem{cor}[theorem]{Corollary}

\newtheorem{rem}{Remark}[section]

 \def\Z{{\mathbb{Z}}}
 \def\Q{{\mathbb{Q}}}
 \def\N{{\mathbb{N}}}

\begin{document}

\title[Homology groups of mapping class groups: a survey]
{Low-dimensional homology groups of mapping class groups: a survey}
\author[KORKMAZ]{ Mustafa Korkmaz}

\address{Department of Mathematics, Middle East Technical University,
 06531 Ankara, Turkey} 
\email{korkmaz@arf.math.metu.edu.tr}

\maketitle
\begin{abstract}
In this survey paper, we give a complete list of known results on the
first and the second homology groups of surface mapping class groups. 
Some known results on higher (co)homology are also mentioned.
\end{abstract}

\section{Introduction}\label{intr}  
 Let $\Sigma_{g,r}^n$ be a connected orientable surface of genus $g$ with $r$
 boundary components and $n$ punctures. The mapping class group of $\Sigma_{g,r}^n$ 
may be defined in different ways. For our purpose, it 
is defined as the group of the isotopy classes of orientation-preserving diffeomorphisms 
$\Sigma_{g,r}^n\to \Sigma_{g,r}^n$. The diffeomorphisms and the isotopies are assumed to
fix each puncture and the points on the boundary. We denote the mapping class group 
of $\Sigma_{g,r}^n$ by $\M_{g,r}^n$.  Here, we see the punctures on the
 surface as distinguished points. 
 If $r$ and/or $n$ is zero, then we omit it from the notation. We write $\Sigma$ for 
the surface $\Sigma_{g,r}^n$ when we do not want to emphasize $g,r,n$.

 The theory of mapping class groups plays a central role in low-dimensional topology.
 When $r=0$ and $2g+n\geq 3$, the mapping class group $\M_g^n$
acts properly discontinuously on the Teichm\"{u}ller space which
 is homeomorphic to some Euclidean space and the stabilizer of each point is finite. The
 quotient of the Teichm\"uller space by the action of the mapping class group is the moduli space
 of complex curves. 

Recent developments in low-dimensional topology made the algebraic structure of the 
mapping class group more important. The examples of such developments are
the theory of Lefschetz fibrations and the Stein fillability of contact $3$-manifolds.
Questions about the structure of these can be stated purely as an algebraic problem
in the mapping class group, but in this paper we do not address such problems.

The purpose of this survey paper is to give a list of complete known results on the homology
groups of the mapping  class groups in dimensions one and two. There is 
no new result in the paper, but there are some new proofs. For example, 
although the first homology group of the mapping class group in genus one case is
known, it seems that it does not appear in the literature.
Another example is that we give another proof of the fact that Dehn twists about nonseparating
simple closed curves are not enough to generate the mapping class group of a surface of genus
one with $r\geq 2$ boundary components (cf. Corollary~\ref{corollaryt} below), as opposed to
the higher genus case.

We shall mainly be interested in orientable surfaces. The first and the second homology groups 
of the mapping class group have been known for more than twenty years. 
We will give the complete list of the first homologies and we calculate them.
An elementary proof of the second homology of the mapping class group was recently
given by the author and Stipsicz in~\cite{ks}. This proof is based on the presentation of
the mapping class group and is sketched in Section~\ref{s5}.
We then outline some known results for higher dimensional (co)homology.
Finally, in the last section, we give the first homology groups of the mapping 
class groups of nonorientable surfaces.

\section{Dehn twists and relations among them}\label{s1} 
Let $\Sigma$ be an oriented surface and let $a$ be a simple closed curve on it.
We always assume that the curves are unoriented.
Cutting the surface $\Sigma$ along $a$, twisting one of the sides by $360$ degrees 
to the right and gluing back gives a self-diffeomorphism of the surface 
$\Sigma$ (cf. Fig.~\ref{sekil1} (a)). Let us denote this
diffeomorphism by $t_a$. In general, a diffeomorphism and its isotopy class
will be denoted by the same letter, so that $t_a$ also represents an element of $\mcg$.
Accordingly, a simple closed curve and its isotopy class are denoted by the same letter.
It can easily be seen that the mapping class $t_a$ depends only on the isotopy class of $a$.
The mapping class $t_a$ is called the (right) Dehn twist about $a$. 

From the definition of a Dehn twist it is easy to see that if
$f:\Sigma\to\Sigma$ is a diffeomorphism and $a$ a simple closed curve on $\Sigma$, then 
there is the equality
\begin{eqnarray}
ft_af^{-1}=t_{f(a)}. \label{eqn=***} 
\end{eqnarray}

 We note that we use the functional notation for the composition
of functions, so that $(fg)(x)=f(g(x))$.

\subsection{The braid relations}
Suppose that $a$ and $b$ are two disjoint simple closed curves on a surface
$\Sigma$. Since the support of the Dehn twist $t_a$ can be chosen to be
disjoint from $b$, we have 
$t_a(b)=b$. Thus by (\ref{eqn=***}), we get
\begin{eqnarray}
t_at_b=t_at_bt_a^{-1}t_a=t_{t_a(b)} t_a=t_bt_a. \label{eqn=commute}
\end{eqnarray}

 \begin{figure}[hb]
 \begin{center}
     \includegraphics[width=13cm]{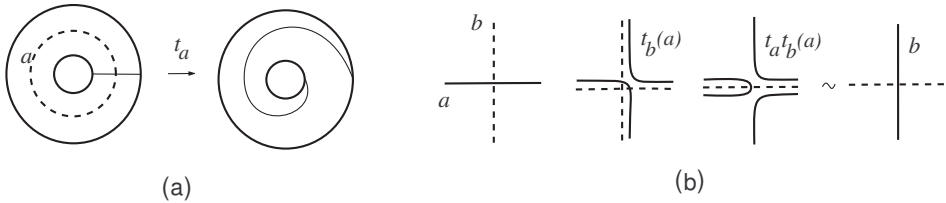}
  \caption{The Dehn twist $t_a$, and the proof of $t_at_b(a)=b$.}
   \label{sekil1}
   \end{center}
 \end{figure}

Suppose that two simple closed curves $a$ and $b$ intersect transversely at only one point.
It can easily be shown that $t_at_b(a)=b$ (cf. Fig.~\ref{sekil1}~(b)). Hence,
\begin{eqnarray}
t_at_bt_a=t_at_bt_a t_b^{-1}t_a^{-1}t_at_b=t_{t_at_b(a)}t_at_b=t_bt_at_b.
\end{eqnarray}

\subsection{The two-holed torus relation}
\begin{figure}[hbt]
 \begin{center}
    \includegraphics{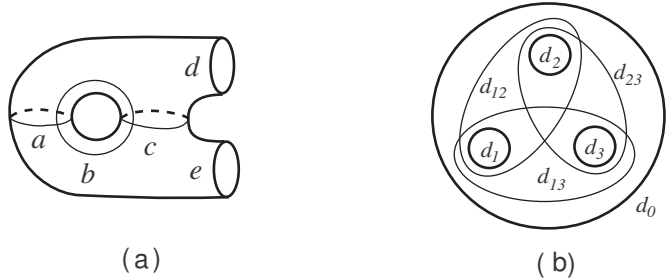}
  \caption{The circles of the two-holed torus and the lantern relations.}
  \label{torus-lantern}
   \end{center}
 \end{figure}

Suppose that $a$, $b$, $c$ are three nonseparating simple closed curves on a 
surface $\Sigma$ such that 
$a$ is disjoint from $c$, and $b$ intersects $a$ and $c$ transversely at one point 
(cf. Fig.~\ref{torus-lantern}~(a)).
A regular neighborhood of $a\cup b\cup c$ is a torus with two boundary components, say
$d$ and $e$. Then the Dehn twists about these simple closed curves satisfy the relation
 \begin{eqnarray}
(t_at_bt_c)^4=t_dt_e. \label{eqn=two-holed}
 \end{eqnarray}
For a proof of this, see~\cite{k1}, Lemma $2.8$. We call it the two-holed torus relation.
In fact, it follows from the braid relations that the three Dehn twists $t_a,t_b,t_c$ on the 
left hand side of this relation can be taken in any order.

 \subsection{The lantern relation}

 Consider a sphere $X$ with four holes embedded in a surface $\Sigma$ in such
 a way that the boundary circles of $X$ are the simple closed curves
 $d_0,d_1,d_2$ and $d_3$ on $\Sigma$ (cf. Fig.~\ref{torus-lantern}~(b)).
 There are three circles $d_{12},d_{13}$ and $d_{23}$
 on $X$ such that there is the relation
 \begin{eqnarray}
 t_{d_0}t_{d_1}t_{d_2}t_{d_3}=t_{d_{12}}t_{d_{13}}t_{d_{23}}. \label{eqn=lantern}
 \end{eqnarray}
This relation is called the lantern relation. It was first discovered by
Dehn~\cite{de} and rediscovered and made popular by Johnson~\cite{j}.

The two-holed torus and the lantern relations can be proved easily: Choose a set 
of arcs dividing the supporting subsurface
into a disc and show that the actions on this set of arcs of the diffeomorphisms on 
the two sides are equal up to homotopy.

 \section{Generating the mapping class group}\label{s2} 
 The search of the algebraic structures of the mapping class group was initiated
by the work of Dehn~\cite{de}. He proved that 
that the mapping class group of a closed orientable surface
is generated by finitely many (Dehn) twists about nonseparating simple closed curves.

\begin{figure}[hbt]
 \begin{center}
    \includegraphics{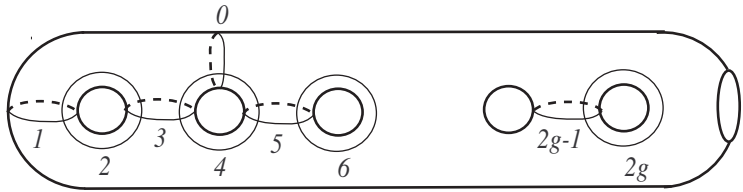}
  \caption{The label $n$ represents the circle $a_n$.}
  \label{generators}
   \end{center}
 \end{figure}

In~\cite{li1,li3}, Lickorish reproduced this result; he proved that 
the mapping class group $\M_g$ can be generated by $3g-1$ Dehn twists,
all of which  are about nonseparating simple closed curves. In~\cite{hu}, Humphries
reduced this number to $2g+1$: The mapping class group $\M_g$ of
a closed orientable surface $\Sigma_g$ of genus $g$ is generated by Dehn twists
about $2g+1$ simple closed curves $a_0,a_1,\ldots,a_{2g}$ of Fig.~\ref{generators}.
In the figure, we glue a disc to the boundary component of the surface
to get the closed surface $\Sigma_g$.
Humphries also showed that the number $2g+1$ is minimal; the mapping 
class group of a closed orientable surface of genus $g\geq 2$ cannot be generated by 
$2g$ (or less) Dehn twists. For any generating set the situation
is different of course; $\M_g$ is generated by two elements. This result is due to Wajnryb~\cite{wa2}.
This is the least number of generators, because $\M_g$ is not commutative.

 Let $\Sigma_g^n$ be an orientable surface of genus $g$ with $n$ punctures.
 Let $n\geq 1$ and let us fix a puncture $x$. By forgetting the puncture $x$,
every diffeomorphism $\Sigma_g^n\to \Sigma_g^n$ induces a diffeomorphism 
$\Sigma_g^{n-1}\to \Sigma_g^{n-1}$.
This gives an epimorphism $\M_g^n\to\M_g^{n-1}$ whose kernel is isomorphic to
the fundamental group of $\Sigma_g^{n-1}$ at the base point $x$ (cf.~\cite{bi1}). 
Therefore, we have a short exact sequence  
\begin{eqnarray}
1\to\pi_1(\Sigma_g^{n-1})\to\M_g^n\to\M_g^{n-1}\to 1. \label{eqn=sequence1}
\end{eqnarray}

Now let $\Sigma_{g,r}^n$ be an orientable surface of genus $g$ with $r$ boundary components
and $n$ punctures.  Assume that $r\geq 1$. Let $P$ be one of the boundary 
components. By gluing a disc $D$ with one puncture along $P$, we get a surface $\Sigma_{g,r-1}^{n+1}$
of genus $g$ with $r-1$ boundary components and $n+1$ punctures.
A diffeomorphism $\Sigma_{g,r}^n\to \Sigma_{g,r}^n$ extends to a 
diffeomorphism $\Sigma_{g,r-1}^{n+1}\to \Sigma_{g,r-1}^{n+1}$
by defining the extension to be the identity on $D$. This way we get an epimorphism
$\M_{g,r}^n\to\M_{g,r-1}^{n+1}$. Note that a Dehn twist on $\Sigma_{g,r}^{n}$ 
along a simple closed curve parallel to $P$ gives a diffeomorphism 
$\Sigma_{g,r-1}^{n+1}\to\Sigma_{g,r-1}^{n+1}$ isotopic to the identity. Essentially, 
this is the only vanishing mapping class under the map $\M_{g,r}^n\to\M_{g,r-1}^{n+1}$. 
More precisely, we have the short exact sequence 
\begin{eqnarray}
1\to\Z \to\M_{g,r}^n\to\M_{g,r-1}^{n+1}\to 1, \label{eqn=sequence2}
\end{eqnarray}
where $\Z$ is the subgroup of $\M_{g,r}^n$ generated by the Dehn twist 
along a simple closed curve parallel to $P$.

It follows from the description of the homomorphisms in the short exact 
sequences~(\ref{eqn=sequence1}) and (\ref{eqn=sequence2}), and the fact that the 
mapping class group $\M_{g}$ is generated by Dehn twists about finitely 
many nonseparating simple closed curves, the group $\M_{g,r}^n$ is generated by 
Dehn twists along finitely many nonseparating simple closed curves and
the Dehn twist along a simple closed curve parallel to each boundary component.

Suppose that the genus of the surface $\Sigma_{g,r}^n$ is at least $2$. The four-holed sphere $X$
of the lantern relation can be embedded in the surface $\Sigma_{g,r}^n$ in such a way that one of the 
boundary components of $X$ is a given boundary component of $\Sigma_{g,r}^n$ and all other six 
curves of the lantern relation are nonseparating on $\Sigma_{g,r}^n$. We conclude from this that

 \begin{thm} \label{thm=generators}
 If $g\geq 2$ then the mapping class group  $ \mcg $ is generated by
Dehn twists about finitely many nonseparating simple closed curves.
\end{thm}

We note that this theorem does not hold  for $g=1$ and $r\geq 2$. 
See Corollary~\ref{corollaryt} in Section~\ref{s3} below.
Dehn twists about boundary  parallel simple closed curves are needed in
this case. In the case that $g=0$, there is no nonseparating simple closed curve.

 \section{Presenting the mapping class groups $\M_{g}$ and $\M_{g,1}$ } \label{s4} 
 
The mapping class groups are finitely presented. The presentation of $\Gamma_2$ 
was first obtained by Birman and Hilden \cite{bh}. For $g\geq 3$, this fact was first proved
by McCool~\cite{mc}
using combinatorial group theory without giving an actual presentation. 
A geometric proof of this was given by Hatcher and 
Thurston~\cite{ht}, again without an explicit presentation. Their proof used the
connectedness and the simple connectedness of a certain complex formed by so called cut systems.
Harer~\cite{ha1} modified the Hatcher-Thurston complex of cut systems in order to
calculate the second homology groups of mapping class groups of orientable surfaces of genus $g\geq 5$.
Using this modified complex, simple 
presentations of the mapping class groups $\M_{g,1}$ and $\M_g$ were finally obtained by Wajnryb~\cite{wa1}. 
Minor errors in~\cite{wa1} were corrected in~\cite{bw}. The proof of Hatcher and Thurston is very 
complicated. In~\cite{wa3}, Wajnryb gave an elementary proof of the presentations of $\M_{g,1}$ and $\M_g$.
This proof does not use the results of Hatcher-Thurston and Harer. It turns out that all the relations 
needed to present the mapping class groups are those given in Section~\ref{s1}, which were obtained 
by Dehn~\cite{de}. 

We now give the Wajnryb presentations of $\M_{g,1}$ and $\M_g$.
So suppose that $n=0$ and $r\leq 1$. As a model for $\Sigma_{g,r}$,
consider the surface in Fig.~\ref{generators}. On the surface $\Sigma_{g,r}$, 
consider the simple closed curves $a_0, a_1,\ldots,a_{2g}$ illustrated in Fig.~\ref{generators}.
 
Let $F$ be the nonabelian free group freely generated by $x_0,x_1,\ldots,x_{2g}$. 
For $x,y\in F$, let $[x,y]$ denote the commutator $xyx^{-1}y^{-1}$.
In the group $F$, we define some words as follows. Let 
$$A_{ij}=[x_i,x_j]$$ 
if the curve $a_i$ is disjoint from the curve $a_j$ in Fig.~\ref{generators}, and let
$$B_0=x_0x_4x_0x_4^{-1}x_0^{-1}x_4^{-1},$$
$$B_i=x_ix_{i+1}x_{i}x_{i+1}^{-1}x_{i}^{-1}x_{i+1}^{-1}$$
for $i=1,2,\ldots, 2g-1$. 
Let us also define the words 
$$C=(x_1x_2x_3)^{4}x_0^{-1}
        (x_4x_3x_2x_1^2x_2x_3x_4)x_0^{-1} (x_4x_3x_2x_1^2x_2x_3x_4)^{-1}$$
and
$$D=x_1 x_3 x_5 w x_0 w^{-1}x_0^{-1} t_2^{-1}x_0^{-1}t_2(t_2t_1)^{-1}x_0 ^{-1}(t_2t_1), $$
where
 $$t_1=x_2 x_1 x_3 x_2, \,\, t_2=x_4 x_3 x_5 x_4,$$
and
 $$w=x_6 x_5 x_4 x_3 x_2 (t_2 x_6 x_5)^{-1} x_0 (t_2 x_6 x_5)
 (x_4x_3x_2x_1)^{-1}.$$

 In the group $F$, we define one more element $E$ to be 
$$ E=[x_{2g+1}, x_{2g}x_{2g-1}\cdots x_3x_2x_1^2x_2
 x_3\cdots x_{2g-1}x_{2g}],$$  
where
 \begin{eqnarray*}
 x_{2g+1}&=& (u_{g-1}u_{g-2}\cdots u_1) x_1  (u_{g-1}u_{g-2}\cdots u_1)^{-1},\\
 u_1&=& (x_1x_2x_3x_4)^{-1}v_1x_4x_3x_2,\\
u_i &=& (x_{2i-1}x_{2i}x_{2i+1}x_{2i+2})^{-1} v_i x_{2i+2}x_{2i+1}x_{2i} \mbox{ for
 $i=2,\ldots,g-1$},  \\
 v_1 &=&  (x_4x_3x_2x_1^2x_2x_3x_4) x_0 (x_4x_3x_2x_1^2x_2x_3x_4)^{-1},\\
 v_i &=& (w_{i}w_{i-1})^{-1} v_{i-1} (w_{i}w_{i-1}) \mbox{ for $i=2,\ldots, g-1$},\\
w_i &=& x_{2i}x_{2i+1}x_{2i-1}x_{2i} \mbox{  for $i=1,2,\ldots, g-1$}.
\end{eqnarray*}

We would like to note that if we define a homomorphism  from $F$ to $\M_{g}$ or  $\M_{g,1}$
by $x_i\mapsto t_{a_i}$, then the relation $C$ maps to a two-holed torus relation and $D$ maps to a 
lantern relation such that all seven simple closed curves in the relation are nonseparating. 
$A_{ij}$ and $B_i$ map to the braid relations.
 
Let us denote by $R_1$ the normal subgroup of $F$  normally generated by the elements $A_{ij}, B_0,B_1,\ldots,
B_{2g-1}, C$ and $D$, and let $R_0$ denote the normal subgroup of $F$ normally generated by
$R_1$ and $E$. The Wajnryb presentation of the mapping class groups $\M_g$ and $\M_{g,1}$
can be summarized as the next theorem.

 \begin{thm}[\cite{wa3}, Theorems $1'$ and $2$] \label{presentation}
 Let $g\geq 2$. Then there are two short exact sequences
 \begin{eqnarray}
1\longrightarrow R_1\longrightarrow F\stackrel{\phi_1}\longrightarrow \M _{g,1}\longrightarrow 1
\label{sequence1}
\end{eqnarray}
and 
 \begin{eqnarray}
1\longrightarrow R_0\longrightarrow F\stackrel{\phi_0}\longrightarrow \M _{g}\longrightarrow 1,
\label{sequence2}
\end{eqnarray}
where $\phi_i(x_j)$ is the Dehn twist $t_{a_j}$ about the curve $a_j$.
 \end{thm}

 \begin{rem}
 If $g=2$ then the relation $D$ is not supported in the surface. In this case, 
one should omit the element $D$ from the definition of $R_0$ and $R_1$.
 \end{rem}

Notice that the presentation of $\M_{2,1}$ in Theorem~$2$ in~\cite{wa3} is slightly
different but equivalent to the presentation above.

A finite presentation of the mapping class group $\M_{g,r}$ is obtained by Gervais in~\cite{ge2}.

 \section{The first homology}\label{s3} 

Recall that for a discrete group $G$, the first homology group $H_1(G;\Z)$ of $G$ with integral coefficients
is isomorphic to the derived quotient $G/[G,G]$, where $[G,G]$ is the subgroup of $G$ generated by 
all commutators $[x,y]$ for $x,y\in G$. Here, $[x,y]=xyx^{-1}y^{-1}$.

From the presentation of the mapping class group $\M_g$ given in Theorem~\ref{presentation}, 
the group  $H_1(\M_g;\Z)$ can be computed easily; it is isomorphic to $\Z_{10}$ if $g=2$ and $0$ if
$g\geq 3$. The fact that $H_1(\M_2;\Z)$ is isomorphic to $\Z_{10}$ was first proved by Mumford~\cite{mu}
and that $H_1(\M_g;\Z)=0$ for $g\geq 3$ by Powell~\cite{po}. We prove this result for $g\geq 3$ 
without appealing to the presentation. We also determine the first homology groups for arbitrary 
$r$ and $n$, which is well known.
 
If $a$ and $b$ are two nonseparating simple closed curves on a surface $\Sigma_{g,r}^n$, then 
by the classification of surfaces there is a
diffeomorphism $f:\Sigma_{g,r}^n\to\Sigma_{g,r}^n$  such that $f(a)=b$. Thus, by~(\ref{eqn=***}), 
we have that $t_b=ft_af^{-1}$. This can also be written as $t_b=[f,t_a]t_a$.
Therefore, $t_a$ and $t_b$ represent the same class $\tau$ 
in $H_1(\M_{g,r}^n;\Z)$.
Since the mapping class group $\M_{g,r}^n$ is generated by Dehn twists about
nonseparating simple closed curves for $g\geq 2$, it follows that the group $H_1(\M_{g,r}^n;\Z)$
is cyclic and is generated by $\tau$. 

\begin{figure}[hbt]
 \begin{center}
    \includegraphics{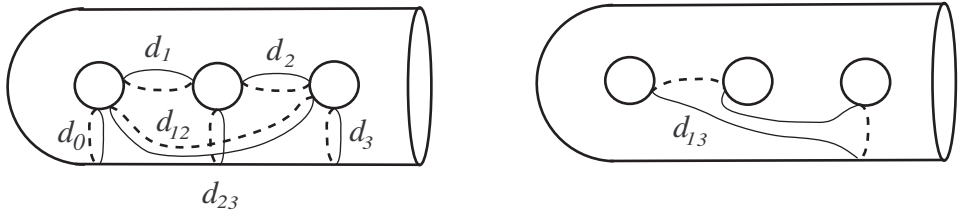}
  \caption{An embedding of the lantern with all curves nonseparating.}
  \label{embed-lantern}
   \end{center}
 \end{figure}

Suppose that $g\geq 3$. The four-holed sphere $X$ of the lantern relation
can be embedded in $\Sigma_{g,r}^n$ such that all seven curves involved in the lantern relation become
nonseparating on $\Sigma_{g,r}^n$ (cf. Fig.~\ref{embed-lantern}). This gives us the relation $4\tau=3\tau$
in $H_1(\mcg;\Z)$. Hence, $H_1(\mcg;\Z)$ is trivial.

Suppose now that $g=2$.  The two-holed torus of the 
two-holed torus relation can be embedded in $\Sigma_{2,r}^n$ so that all five curves in the relation
becomes nonseparating on $\Sigma_{2,r}^n$. This gives us $12\tau=2\tau$, i.e. $10\tau=0$.
On the other hand, since there is an epimorphism $\M_{2,r}^n\to  \M_2$, it follows that
$H_1(\M_{2,r}^n)=\Z_{10}$.

Although we usually deal with the surfaces of genus at least two, we would like to mention
the first homology groups in the genus one case as well. 

Consider a torus $\Sigma_{1,r}^n$ with $n$ punctures and $r$ boundary components, 
$P_1,P_2,\ldots,P_r$. For each $i=1,2,\ldots, r$, let $\partial_i$ be a simple closed curve 
parallel to $P_i$.  The mapping class group $\M_{1}$ is 
generated by the Dehn twists about two (automatically nonseparating) simple closed curves intersecting
transversely at one point. It can be proved by the exact sequences 
(\ref{eqn=sequence1}) and  (\ref{eqn=sequence2}) that the mapping class group 
$\Sigma_{1,r}^n$ is generated by Dehn twists about finitely many nonseparating simple closed curves 
and $r$ Dehn twists about  $\partial_1,\partial_2,\ldots,\partial_r$.
By the use of the lantern relation, it can be shown that we may omit any one of $\partial_i$, say $\partial_n$.
It will follow from the following computation of the first homology of $\M_{1,r}^n$ that in fact 
no more $\partial_i$ may be omitted. 

By the argument in the case of higher genus, any two Dehn twists about nonseparating simple 
closed curves are conjugate. Hence, they represent the same class $\tau$ in 
 $ H_1(\M_{1,r}^n;\Z)$.

Assume first that $r=0$. The group  $\M_1$ is isomorphic to
$SL(2,\Z)$. Hence, $H_1(\M_1;\Z)$ is isomorphic to $\Z_{12}$ and is generated by $\tau$.
Since  $\M_1^n$ is generated by Dehn twists about nonseparating simple closed
curves, the homology group $H_1(\M_1^n;\Z)$ is cyclic and generated by $\tau$.
It was shown in Theorem~$3.4$ in~\cite{km} that $12\tau=0$. On the other hand, 
the surjective homomorphism $\M_1^n\to \M_1$ obtained by forgetting the punctures
induces a surjective homomorphism between  the first 
homology groups, mapping $\tau$ to the generator of $H_1(\M_1;\Z)$.
It follows that $H_1(\M_1^n;\Z)$ is isomorphic to $\Z_{12}$.

The group $\M_1^1$ is also isomorphic to $SL(2,\Z)$.
Let $a$ and $b$ be two simple closed curve on $\Sigma_{1,1}$ intersecting each other
transversely at one point. By examining the short exact sequence 
$$
1\to\Z \to \M_{1,1} \to \M_1^1\to 1 
$$ 
it can be shown easily that $\M_{1,1}$ has a presentation with generator 
$t_a,t_b$ and with a unique relation $t_at_bt_a=t_bt_at_b$. That is, $\M_{1,1}$ is isomorphic to
the braid group on three strings. Hence, by abelianizing this presentation,
we see that  $H_1(\M_{1,1};\Z)$ is isomorphic to $\Z$.

Assume now that $r\geq 1$. Since $\M_{1,r}^n$ is generated by Dehn twists
about nonseparating simple closed curves and the curves
$\partial_1,\partial_2,\ldots,\partial_{r-1}$, the group
$H_1(\M_{1,r}^n;\Z)$ is generated by $\tau,\delta_1,\delta_2,\ldots,\delta_{r-1}$, 
where $\delta_i$ is the class in $H_1(\M_{1,r}^n;\Z)$ of the Dehn twist about $\partial_i$.

We prove that $\tau,\delta_1,\delta_2,\ldots,\delta_{r-1}$ are linearly independent.
Let $n_0 \tau +n_1 \delta_1 +n_2 \delta_2 +\ldots +n_{r-1}\delta_{r-1}=0$ with $n_i\in\Z$. Gluing a 
disc to each $P_i$ for $i\leq r-1$ and forgetting the punctures gives rise to an epimorphism  
$H_1(\M_{1,r}^n;\Z) \to H_1(\M_{1,1};\Z)$.  Under this map, 
$\tau$ is mapped to the generator of $H_1(\M_{1,1};\Z)$
and all $\delta_i$ to zero. Hence, $n_0=0$. 
Similarly, gluing a  disc to each boundary component but $P_i$ and forgetting the punctures
induces an epimorphism $H_1(\M_{1,r}^n;\Z) \to H_1(\M_{1,1};\Z)$ mapping
$\delta_i$ to $12\tau$ and each $\delta_j$, $j\neq i$, to $0$, where $\tau$ is 
the generator of $H_1(\M_{1,1};\Z)$. This shows that $n_i=0$ for each $i=1,2,\ldots, r-1$.
We conclude that $H_1(\M_{1,r}^n;\Z)$ is isomorphic to $\Z^r$. 

We collect the results of this section in the next theorem. 

\begin{thm}
Let $g\geq 1$. The first homology group $H_1(\M_{g,r}^n;\Z)$ of the mapping 
class group is isomorphic to $\Z_{12}$ if $(g,r)=(1,0)$,  $\Z^r$ if $g=1,r\geq 1$,
$\Z_{10}$ if $g=2$ and $0$ if $g\geq 3$.
\end{thm}

\begin{cor} \label{corollaryt}
Let $r\geq 2$. The mapping class group $\M_{1,r}^n$ cannot
be generated by Dehn twists about nonseparating simple closed curves.
\end{cor}

We note that this corollary was proved by Gervais in~\cite{ge1} by a different argument. 

\bigskip

 \section{The second homology}\label{s5} 

The second homology group of the mapping class group $\M_{g,r}^n$ 
for $g\geq 5$ was first computed by Harer in~\cite{ha1}. His proof relies
on the simple connectedness of a complex obtained by modifying
the Hatcher-Thurston complex~\cite{ht}.  But
this proof is extremely complicated to understand. The computation of
$H_2(\M_g;\Z)$ in~\cite{ha1} was incorrect and it was corrected later. See,
for examle, \cite{ha2} or \cite{mo1}.

In~\cite{pi}, Pitsch gave a simple proof of $H_2(\M_{g,1};\Z)=\Z$ for $g\geq 4$.
His method used the presentation of the mapping class group $\M_{g,1}$
and the following theorem of Hopf (cf.~\cite{br}): Given a short exact sequence of groups
$$1\to R\to F\to G\to 1,$$
where $F$ is free, then 
\begin{eqnarray}
H_2(G;\Z) =\frac{R\cap[F,F]}{[R,F]}.
\label{H_2}
\end{eqnarray}

In~\cite{ks}, Stipsicz and the author extended Pitsch's proof to 
$H_2(\M_{g};\Z)=\Z$ for $g\geq 4$. Then the homology 
stabilization theorem of Harer~\cite{ha2}  and a use of the Hochschild-Serre
spectral sequence for group extensions
enabled us to give a new proof of Harer's theorem
on the second homology of mapping class groups, by extending it to the $g=4$ case.

We now outline the proof of $H_2(\M_{g};\Z)=\Z$ for $g\geq 4$.

Consider the short exact sequence~(\ref{sequence2}). Recall that $F$ in ~(\ref{sequence2}) 
is the free group generated freely by $x_0,x_1,\ldots, x_{2g}$ and $R_0$ is the
normal subgroup of $F$ normally generated by the elements $A_{ij}, B_0,B_1,\ldots,B_{2g-1},C,D$ and $E$.
By Hopf's theorem, we have $$ H_2(\M_g;\Z) =\frac{R_0\cap[F,F]}{[R_0,F]}.$$ 
Hence, every element in $H_2(\M_g;\Z)$ has a representative of the form 
\begin{eqnarray}
AB_0^{n_0}B_1^{n_1}\cdots B_{2g-1}^{n_{2g-1}} C^{n_C} D^{n_D}E^{n_E}, \label{eqn=ABCDE}
 \end{eqnarray}
where $A$ is a product of $A_{ij}$. 

Note that each $A_{ij}$ and $E$ represent elements of $H_2(\M_g;\Z)$ since they are contained in
$R_0\cap [F,F]$.

It was shown in~\cite{pi} by the use of the lantern relation that each $A_{ij}$ represents
the trivial class in $\displaystyle H_2(\M_{g,1};\Z)=(R_1\cap [F,F])/[R_1,F]$. 
The same proof applies to show that the class of each $A_{ij}$ in $H_2(\M_g;\Z)$ is zero.
The main reason for this is that since $g\geq 4$, for any nonseparating
simple closed curve $a$ on $\Sigma_g$, the four-holed sphere of the lantern relation
can be embedded in $\Sigma_g - a$ such that all seven curves of the relation are nonseparating
on $\Sigma_g-a$.

The main improvement in~\cite{ks} after~\cite{pi} is to show that $E$ represents 
the zero element in $H_2(\M_g;\Z)$. The proof of this uses the braid relations and 
the two-holed torus relation. Therefore, we may delete $A$ and $E$ in~(\ref{eqn=ABCDE}).

Note that an element of $F$ is in the derived subgroup $[F,F]$ of $F$
if and only if the sum of the exponents of each generator $x_i$ is zero.
Since the expression~(\ref{eqn=ABCDE}) must be in $[F,F]$, by looking at the sum of 
the exponents of the generators $x_{2g},x_{2g-1},\ldots, x_{6}$, one can see 
easily that $n_{2g-1},n_{2g-2},\ldots,n_5$ must be zero.
Then, by looking at the sums of the exponents of the other generators,  
it can be concluded that there must be 
an integer $k$ such that $n_0=-18k, n_1=6k,n_2=2k,n_3=8k,n_4=-10k,n_C=k$ and $n_D=-10k$.
This says that $H_2(\M_{g};\Z)$ is cyclic and is generated
by  the class of the element
\begin{eqnarray}
 B_0^{-18}B_1^{6}B_2^{2} B_3^{8}B_4^{-10}CD^{-10}. \label{eqn=BBB}
 \end{eqnarray}

On the other hand, for every $g\geq 3$, the existence of a genus-$g$ surface bundle 
with nonzero signature guarantees that $H^2(\M_g;\Z)$ contains an infinite cyclic subgroup
(cf.~\cite{me}); the signature cocycle is of infinite order.
The universal coefficient theorem implies that 
$H_2(\M_g;\Z)$ contains an element of infinite order. This shows that 
$H_2(\M_g;\Z)=\Z$ for $g\geq 4$.

By omitting $E$ from the above proof, the same argument also proves that 
$H_2(\M_{g,1};\Z)=\Z$ for $g\geq 4$.

A special case of Harer's  homology stability theorem in~\cite{ha2} says that 
for $g\geq 4$ and $r\geq 1$ the inclusion mapping $\Sigma_{g,r}\to\Sigma_{g,r+1}$ obtained by 
gluing a disc with two boundary components to $\Sigma_{g,r}$ along one of the boundary
components induces an isomorphism $H_2 (\M_{g,r};\Z )\to H_2(\M_{g,r+1};\Z).$ 
Also, an application of the Hochschild-Serre spectral sequence to~(\ref{eqn=sequence2})
shows that $H_2 (\M_{g,r}^n ;\Z )=H_2 (\M_{g,r+1}^{n-1};\Z )\oplus \Z$
for $g\geq 3$.

We can summarize the results mentioned above as follows. The details of the proof may
be found in~\cite{ks}.

\begin{thm}
If $g\geq 4$ then $H_2(\M_{g,r}^n;\Z)$ is isomorphic to $\Z^{n+1}$.
\end{thm}

The same method above also proves that $H_2(\M_2;\Z)=H_2(\M_{2,1};\Z)$ is isomorphic to
either $0$ or $\Z_2$, and the groups $H_2(\M_3;\Z)$ and $H_2(\M_{3,1};\Z)$ are isomorphic to
either $\Z$ or $\Z\oplus \Z_2$. By the work of Benson-Cohen~\cite{bc}, 
$H_2(\M_2;\Z_2)$ is isomorphic to $\Z_2\oplus\Z_2$. It follows now from the universal coefficient
theorem that $H_2(\M_2;\Z)$ is not trivial, hence $\Z_2$.  To the best knowledge of the
author, the computations of $H_2(\M_{2,r}^n;\Z)$ in the remaining cases and 
$H_2(\M_{3,r}^n;\Z)$ are still open.

 \section{Higher (co)homologies}\label{s6} 

Here we will mention a few known results on the (co)homology group of the mapping class group. 
In Section~\ref{s5}, we appealed to a special case of the homology stability theorem of Harer.
The original theorem asserts that in a given dimension the homology group of the mapping 
class group of a surface of with boundary components does not depend on the genus if 
the genus of the surface is sufficiently high. 
This result was improved by Ivanov in~\cite{iv1,iv2}. In~\cite{iv2},
Ivanov also proved a stabilization theorem for the homology with twisted coefficients
of the mapping class groups of closed surfaces.
 
The third homology group of $\M_{g,r}$ with rational coefficients was computed by
Harer in~\cite{ha3}. It turns out that $H_3(\M_{g,r};\Q)=0$ for $g\geq 6$. 

Let $\Q [z_2,z_4,z_6,\ldots]$ denote the polynomial algebra of generators
$z_{2n}$ in dimensions $2n$ for each positive integer $n$. Then there are classes
$y_2,y_4,y_6,\ldots$ with $y_{2n}\in H^{2n}(\M_g;\Q )$ such that the homomorphism 
of algebras 
$$\Q [z_2,z_4,z_6,\ldots]\to H^* (\M_g;\Q ) $$
given by $z_{2n}\mapsto y_{2n}$ is an injection in dimensions less than $g/3$.
This result was proved by Miller~\cite{mi}.

The entire mod-$2$ cohomology of $\M_2$ is also known. Benson and Cohen~\cite{bc} computed
the Poincar\'e series for mod-$2$ cohomology to be
$$ (1+t^2+2t^3+t^4+t^5)/(1-t)(1-t^4)=1+t+2t^2+4t^3+6t^4+7t^5+\cdots .$$

\section{Nonorientable surfaces}

In this last section, we outline the known results about the generators and the homology groups
of the mapping class groups of nonorientable surfaces. 
So let $S_g^n$ denote a nonorientable surface of genus $g$ with $n$ punctures.
Recall that the genus of a closed nonorientable surface is defined as the number of real projective planes in a 
connected sum decomposition. Let us define the mapping class group $\M_g^n$ 
as in the orientable case; diffeomorphisms and isotopies are required to fix each puncture.

The mapping class group $\M_1$ of the real projective plane is trivial and the group
$\M_2$ is isomorphic to $\Z_2\oplus\Z_2$ (cf.~\cite{li2}). Lickorish~\cite{li2,li4} and Chillingworth~\cite{ch}
proved that 
if $g\geq 3$ then $\M_g$ is generated by a finite set consisting of Dehn twists about two-sided nonseparating
simple closed curves and a crosscap slide (or Y-homeomorphism). See also~\cite{bich}.
Using this result the author~\cite{k1}
computed $H_1(\M_g;\Z)$. This result was extended to the punctured cases in \cite{k3}. We note that
the group $\M_g^n$ of this section is called the pure mapping class group in~\cite{k3} and
denoted by ${{\mathcal{PM}}}_{g,n}$. The first homology group of $\M_g^n$ with integer coefficients  is 
as follows.

\begin{thm} [\cite{k3}]
Let $g\geq 7$. Then the first homology group $H_1(\M_g^n;\Z)$ of $\M_g^n$ is isomorphic to $\Z_2^{n+1}$.
\end{thm}

If we define ${\mathcal M}_g^n$ as the group of the diffeomorphisms $S_g^n\to S_g^n$
modulo the diffeomorphisms which are isotopic to the identity by an isotopy fixing 
each puncture, then more is known. Let $\N_+$ and $\N$ denote the set of positive integers
and the set of nonnegative integers, respectively. Define a function $k:\N_+\times \N\to \N$
by declaring
$k(1,0)=0$, $k(4,0)=3$, $k(g,0)=2$ if $g=2,3,5,6$,
$k(g,0)=1$ if $g\geq 7$, $k(g,1)=k(g,0)+1$ and $k(g,n)=k(g,0)+2$ if $n\geq 2$.

\begin{thm} [\cite{k3}]
The first homology group $H_1({\mathcal M}_g^n;\Z)$ of the mapping class group ${\mathcal M}_g^n$
of a nonorientable surface of genus $g$ with $n$ punctures 
is isomorphic to the direct sum of $k(g,n)$ copies of $\Z_2$.
 \end{thm}

It is easy to see that the groups $\M_g^n$ and ${\mathcal M}_g^n$ fit into a short exact sequence
$$
1\longrightarrow \M_g^n \longrightarrow {\mathcal M}_g^n
\longrightarrow Sym(n) \longrightarrow 1,
$$
where $Sym(n)$ is the symmetric group on $n$ letters.

No higher homology groups of $\M_g^n$ or ${\mathcal M}_g^n$
are known.

\bigskip

{\bf Acknowledgements:} The author wishes to thank the referee for his/her suggestions
on the earlier version of the paper and Selman Akbulut for his encouragement to write such a survey paper.

\end{document}